\UseRawInputEncoding
\documentclass[12pt,thmsa,emstex]{article}
\usepackage{amsfonts}
\usepackage{t1enc}
\usepackage{amsmath}
\usepackage{amssymb}
\usepackage{epsfig}
\usepackage[utf8]{inputenc}							
\usepackage[french,english]{babel}
\usepackage{mathrsfs}
\usepackage{hyperref}
\usepackage{amsthm}
\usepackage{graphics}
\usepackage{layout}
\usepackage{latexsym}
\usepackage[dvipsnames]{xcolor}
\usepackage{indentfirst} 						
\usepackage{yfonts}

\usepackage[a4paper]{geometry}
\def\mylabel#1{\label{#1}}

\newtheorem{theorem}{Theorem}[section]
\newtheorem{lemma}[theorem]{Lemma}
\newtheorem{corollary}[theorem]{Corollary}
\newtheorem{proposition}[theorem]{Proposition}
\newtheorem{example}[theorem]{Example}
\newtheorem{examples}[theorem]{Examples}
\newtheorem{remark}[theorem]{Remark}

\newtheorem{definition}[theorem]{Definition}
\newtheorem{exercice}[theorem]{Exercice}

\DeclareMathOperator{\dist}{dist}

\def\bit{\begin{itemize}}

\def\eit{\end{itemize}}

\reversemarginpar   

\def\bc{\begin{center}}

\def\ec{\end{center}}

\def\bthm{\begin{theorem}}

\def\ethm{\end{theorem}}

\def\bcor{\begin{corollary}}

\def\ecor{\end{corollary}}

\def\bprop{\begin{proposition}}

\def\eprop{\end{proposition}}

\def\blem{\begin{lemma}}

\def\elem{\end{lemma}}

\def\bex{\begin{example}}

\def\eex{\end{example}}

\def\bexs{\begin{examples}}

\def\eexs{\end{examples}}

\def\bexo{\begin{exercice} \rm }

\def\eexo{\end{exercice} }

\def\brem{\begin{remark}}

\def\erem{\end{remark}}

\def\prf{{\bf Proof }}

\def\bdes{\begin{description}}

\def\edes{\end{description}}

\def\ita{\item[(a)]}

\def\itb{\item[(b)]}

\def\itc{\item[(c)]}

\def\beq{\begin{equation}}

\def\eeq{\end{equation}}

\def\ben{\begin{enumerate}}

\def\een{\end{enumerate}}

\def\beqar{\begin{eqnarray}}

\def\eeqar{\end{eqnarray}}

\def\beqarr{\begin{eqnarray*}}

\def\eeqarr{\end{eqnarray*}}

\def\prf{{\bf Proof }\hspace{.1in}}

\def\PP{{\mathcal P}}
\def\E{{\mathbb E}}

\def\P{{\mathbb P}}

\def\Ind{{\mathbf 1}}
\def\RR{{\mathbb R}}  

\def\CCC{{\mathcal C}}
\def\NN{{\mathbb N}}

\def\eps{\varepsilon}

\begin{document}
\title{Stochastic persistence in degenerate stochastic Lotka-Volterra food chains}
\author{Michel Benaïm $^1$, Antoine Bourquin $^1$, Dang H. Nguyen $^2$
\bigskip
\\ \normalsize $^1$ Institut de Math\'ematiques, Universit\'e de Neuch\^atel, Switzerland
\\ \normalsize $^2$ Department of Mathematics, University of Alabama, USA}


\maketitle

\begin{abstract}
We consider a Lotka-Volterra food chain model with possibly intra-specific competition in a stochastic environment represented by stochastic differential equations. In the non-degenerate setting, this model has already been studied by A. Hening and D. Nguyen in \cite{HeningDang18c,HeningDang18b} where they provided conditions for stochastic persistence and extinction.
In this paper, we extend their results to the degenerate situation in which the top or the bottom species is subject to random perturbations. Under the persistence condition, there exists a unique invariant probability measure supported by the interior of $\RR_+^n$ having a smooth density.

Moreover, we study a more general model, in which we give new conditions which make it possible to characterize the convergence of the semi-group towards the unique invariant probability measure either at an exponential rate or at a polynomial one. This will be used in the stochastic Lotka-Volterra food chain to see that if intra-specific competition occurs for all species, the rate of convergence is exponential while in the other cases it is polynomial.
\end{abstract}



\tableofcontents



\section{Introduction} \mylabel{sec:intro}

We consider a stochastic food chain model with $n$ species.  Species $n$ is the apex predator which hunts species $n-1$ which is in turn the predator of species $n-2$ which hunts species $n-3$ and so on until species $1$ which is at the bottom of the food chain.

More precisely, the model is given by the following (possibly degenerate) Lotka-Volterra food chain stochastic differential equation
on $\RR_+^n := \left\{x \in \RR^n \mid x_i \geq 0, \; i = 1,\ldots,n \right\}$,
\begin{align} \label{model_normal}
d X_i(t) = X_i(t) F_i(X(t))\, dt + \sigma_i  X_i(t)\, dB_t^i \qquad i = 1, \ldots, n
\end{align}
whereby
\begin{align} \label{fucntions_F_i}
F_i(x) = \left\{
     \begin{array}{ll}
      a_{10} -a_{11} x_1 - a_{12} x_2 								& i=1,\\
     -a_{i0} + a_{i,i-1}x_{i-1} - a_{ii}\, x_i - a_{i,i+1}x_{i+1} 	& i= 2, \ldots,n-1, \\
      -a_{n0} + a_{n,n-1} x_{n-1} -a_{nn} \,x_n 					& i=n,
     \end{array}
     \right.
\end{align}

with $a_{11},a_{ij} >0$ for $i \neq j$, $a_{ii} \geq 0$, $\sigma_i \geq 0$ for $i=1, \ldots, n$ and $(B_t^1, \ldots, B_t^n)_{t \geq 0}$ a vector of independent Brownian motions.

The quantity $X_i(t)$ represents the density of species $i$ at time $t$, $a_{ii}$ its intra-specific competition rate, $a_{i0}$ its death rate (or growth rate for species $1$), $a_{i,i-1}$ the rate at which it hunts species $i-1$, and $a_{i,i+1}$ the rate at which it is hunted by species $i+1$.
In addition, the capital growth rate of species $i$ is subject to environmental fluctuations whose variance is measured by $\sigma_i^2$.

In the deterministic case (i.e. $\sigma_i=0$ for all $i$), this model was first considered by T. Gard and T. Hallam in the late seventies in \cite{gard1979persistence} when there is no intra-specific competition (i.e $a_{ii}=0$ for all $i\geq 2$). They gave a criterion based uniquely on the coefficients $a_{ij}$ for persistence and extinction. 

Equation \eqref{model_normal} with non-zero $\sigma_i$s, is a natural way to incorporate random fluctuations in the growth rate.
%
Although the choice of a constant variance
$\sigma_i^2$
seems oversimplified as commented on in \cite{allen2016environmental}, stochastic models with linear diffusion terms have been used extensively by biologists for predicting
population dynamics and estimating extinction likelihoods in \cite{dennis1991estimation, foley1994predicting, lande2003stochastic}. They were also used to explore
how stochasticity can facilitate or inhibit the persistence of populations, coexistence of interacting species or genotypes; see
e.g. \cite{evans2013stochastic,schreiber2012evolution} and the references therein.

Following this line, A. Hening and D. Nguyen in \cite{HeningDang18c,HeningDang18b} studied \eqref{model_normal} in the non-degenerate case (i.e. $\sigma_i \neq 0$ for all $i$) and provided general conditions for persistence and extinction.

The main purpose of this paper is to extend their results to the degenerate setting whereby some but not all of the $\sigma_i$s are $0$.
This is motivated by the fact that environmental fluctuations may have different impacts to different species in an ecological system. For instance, some species, perhaps because of their biological nature, may be negligibly affected by the sources of environmental randomness under consideration. 

We will focus our attention here on the specific case whereby the noise affects at least the top or bottom species (meaning that $\sigma_1 \neq 0$ or $\sigma_n \neq 0$) and show that the conditions given in \cite{HeningDang18c,HeningDang18b} ensuring persistence of the process remain valid. Here by persistence we mean that:
 \begin{description}
 \item[(a)]
 There exists a unique invariant probability measure on $\RR_{++}^n := \{ x \in \RR_+^n \mid \prod_{i=1}^n x_i >0 \}$ called the {\em persistent measure}.
 \item[(b)] Whenever the initial condition lies in $\RR_{++}^n,$ the process converges (in distribution) to the persistent measure.
 \end{description}

 Another key question is to understand the rate of this convergence and how this rate is determined by the model parameters.
 We will see here that  positive intra-specific rates (i.e $a_{ii} > 0$ for all $i$)  ensure an exponential rate of convergence (exactly like in the non-degenerate case (\cite{HeningDang18} \cite{HeningDang18c}))  whilst, if some of the rates are zero (i.e~ $a_{ii} = 0$ for some $i$), one can only ensure a polynomial rate of convergence. This will be established in Section 4, in a general framework providing new results that can be applied  beyond  Lotka-Volterra food chain models.

\paragraph*{Outline of contents} The first part of Section 2 is devoted to some preliminaries and notations. We present the main results of this paper in the second part. More precisely, under some conditions, we characterize the persistence of all species and the rates of convergence to a unique invariant probability measure. We also provide criteria for the extinction of some species. The first part of Section 3 is devoted to the introduction of some mathematical tools and follows \cite{B2018}. In particular, we provide a very useful criterion which ensures the persistence. In the second part, we prove Theorem \ref{theorem_stoch_pers}, which provides conditions for stochastic persistence.

In section 4,  we shall prove general results on the rates of convergence of the semi-group towards the unique invariant probability measure. We will then use these to prove Theorem \ref{thm_rate_of_convergence}, which characterizes the convergence rates of our model. We will see that the intra-specific competition will give an exponential rate whereas without it we will only have a polynomial one. Section 5 then focuses on the case of extinction of at least one species. The appendix (Section 6) contains some technical proofs of Section 2 and 3 for convenience.


\section{Notations and results} \mylabel{main_result}

\subsection{Notations}

Throughout the paper, we let
$$\RR_{++}^n := \{ x \in \RR_+^n \mid \prod_{i=1}^n x_i >0 \} \quad\text{ and }\quad \partial \RR_
+^n := \{ x \in \RR_+^n \mid \prod_{i=1}^n x_i = 0 \}$$
denote, respectively, the interior and the boundary of $\RR_+^n$. For $1\leq j<n$, we define $\RR_+^j$ and $\RR_{++}^j$ in the same way and we also define 
$$\RR_+^{(j)} := \{ x \in \RR_+^n \mid x_1,\ldots,x_j\geq 0,\, x_{j+1},\ldots,x_n=0 \} $$
and
$$ \RR_{++}^{(j)} := \{ x \in \RR_+^n \mid x_1,\ldots,x_j> 0,\, x_{j+1},\ldots,x_n=0 \}.$$

We let $(X^x(t))_{t\geq 0}$ denote the solution to \eqref{model_normal} with initial condition $X^x(0) = x$
and $(P_t)_{t \ge 0}$ denote its transition kernel defined by $P_t f(x) := \E [f(X^x(t))]$ for every measurable bounded function $f : \RR_+^n \to \RR$ and by $P_t(x,B) := P_t \Ind_B (x)$ for all Borel set $B \subset \RR_+^n$.

A Borel set $B$ is called {\it invariant} if $$P_t \Ind_B = \Ind_B$$ for all $t \geq 0$.

Given an invariant set $A \subset \RR_+^n$ (typically $A= \RR_{++}^n$ or $A = \partial \RR_+^n$), we let $\PP_{inv}(A)$ (resp. $\PP_{erg}(A)$) denote the {\it set of invariant} (resp. {\it ergodic}) {\it probability measure} of the process that are supported by $A$, i.e. $\mu(A)=1$ for $\mu \in \PP_{inv}(A)$ (resp. for $\mu \in \PP_{erg}(A)$). Recall that an invariant probability measure $\mu$ is {\it ergodic} if  for every invariant Borel set $B$, $\mu(B)\in\{0,1\}$.

We also let $\left(\Pi_t^x\right)_{t \in \RR_+}$ denote the {\it set of empirical occupation measures} of the process $(X^x(t))$ with initial condition $x \in \RR_+^n$, that is
$$\Pi_t^x (\cdot) := \frac{1}{t} \int_0^t \Ind_{\{ X^x(s) \in \,\cdot \,\}} ds.$$

For a function $f:\RR_+^n \to \RR$ and a measure $\mu$, we write $\mu f := \int f(x) \mu(dx)$.

Recall that the {\it total variation distance} between two probability measures $\mu,\nu $ on $\RR_+^n$ is defined by
$$\|\mu - \nu \|_{TV} := \sup\{|\mu f - \nu f| \mid f:\RR_+^n \to \RR \text{ measurable bounded}, \|f\|_\infty <1\}.$$
In the same manner, for $f:\RR_+^n \to \RR_+$, we define the {\it $f$-norm} by
$$\|\mu \|_f := \sup\{|\mu g |\mid g:\RR_+^n \to \RR \text{ measurable bounded}, |g| <f\}.$$
Remark that when $f=1$, it is the total variation norm.

We say that a map $f : \RR_+^n \to \RR$ is {\it proper} if $\lim\limits_{\|x \| \to \infty} f(x) =  \infty$.

\subsection{Preliminaries}


We start this part by rewriting \eqref{model_normal} in the Stratonovich formalism
\begin{align} \label{model_stratonovich}
dX_i(t) = X_i(t) \widetilde{F}_i (X(t))\, dt +  \sigma_i X_i(t) \circ dB_t^i \qquad i = 1, \ldots, n
\end{align}
where
\begin{align} \label{F_tilde}
\widetilde{F}_i(x)  = F_i (x)-  \frac{1}{2} \sigma_i^2 \qquad i = 1, \ldots, n.
\end{align}

\begin{remark}
If we
let $\widetilde{a}_{10} = a_{10} - \frac{1}{2} \sigma_1^2$, $\widetilde{a}_{i0} = a_{i0} + \frac{1}{2} \sigma_i^2$ for $i \geq 2$ and $\widetilde{a}_{ij} = a_{ij}$ if $i\neq 0$, then $\widetilde{F}$ is defined by \eqref{fucntions_F_i} with $a_{ij}$ replaced by $\widetilde{a}_{ij}$.
\end{remark}

An equilibrium $p^*$ for $\widetilde{F}$ given by \eqref{F_tilde} (i.e. $\widetilde{F}(p^*) = 0$) is called {\it positive} if $p^* \in \RR_{++}^n$. It turns out that the existence of a such point (which has to be unique) is a necessary and sufficient condition for stochastic persistence. Before stating our main theorems, we will give another equivalent  condition of the existence of this positive equilibrium for $\widetilde{F}$ and we start by defining a parameter, called $\widetilde{\delta}(n)$.


For $a,b, m \in \NN$, $a \leq b$, we let $A_a^b (m)$ denote the set of permutations of $ \{a,a+1, \ldots, b-1, b\}$ that are product of $m$ disjoint transpositions of the form  $(i\; i+1)$ and we let
$$A_a^b := \bigcup_{m=0}^\infty \, A_a^b (m).$$
Remark that $A_a^b(m) = \emptyset$ if $m > \frac{b-a+1}{2}$. Now, we define
$$\widetilde{\delta}(n) := \widetilde{a}_{10} \prod_{j=2}^n \widetilde{a}_{j,j-1} - \sum_{k=2}^n \widetilde{a}_{k0} \prod_{l=k+1}^n \widetilde{a}_{l,l-1} \sum_{\alpha \in A_1^{k-1}} \prod_{j=1}^{k-1} \widetilde{a}_{j,\alpha(j)}.$$

\bexs
\begin{align*}
A_1^3 &= \left\{ Id_{\{1,2,3\}}, (1 \; 2), (2 \;3) \right\},  \\
A_1^4 &= \left\{ Id_{\{1,2,3,4\}}, (1 \; 2), (2\; 3), (3 \; 4), (1 \; 2)(3 \; 4) \right\}, \\
\widetilde{\delta}(2) &= \widetilde{a}_{10}\, \widetilde{a}_{21} - \widetilde{a}_{20}\, \widetilde{a}_{11},\\
\widetilde{\delta}(3) &= \widetilde{a}_{10} \,\widetilde{a}_{32} \,\widetilde{a}_{21} - \widetilde{a}_{20} \,\widetilde{a}_{32} \,\widetilde{a}_{11} - \widetilde{a}_{30} \,\widetilde{a}_{22} \,\widetilde{a}_{11} - \widetilde{a}_{30} \,\widetilde{a}_{12} \,\widetilde{a}_{21}.
\end{align*}
\eexs

\brem \label{delta_cst_kappa}
If $\widetilde{a}_{ii} = 0$ for $i \ge 2$, then $\widetilde{\delta}(n) = \widetilde{\kappa}(n) \prod\limits_{i=2}^n \widetilde{a}_{i,i-1}$ where $\widetilde{\kappa}(n)$ is introduced in \cite{HeningDang18b}. In some sense, $\widetilde{\delta}(n)$ can be seen as a generalisation of $\widetilde{\kappa}(n)$, since it is also defined for $a_{ii} \geq 0$.
\erem

We can now state the link between the constant $\widetilde{\delta}(n)$ and the nature of the equilibrium of $\widetilde{F}$.


\begin{proposition} \label{Prop_delt_pos_iff_F_positive_Zero}
$\widetilde{\delta}(n) >0 \Longleftrightarrow \widetilde{F}$ has a positive (necessarily unique) equilibrium $p^*$.
\end{proposition}

The proof is technical and postponed to Appendix \ref{sec_proof_delta_pos}. Moreover, since the existence of a positive equilibrium for $\widetilde{F}$ is a necessary and sufficient condition for stochastic persistence, so is the positivity of $\widetilde{\delta}(n)$.

\begin{remark}
The proof of Proposition \ref{Prop_delt_pos_iff_F_positive_Zero} provides an explicit formula for $p^*$. This extends the results in \cite{HeningDang18c} since it also provides an explicit formula for $\mathcal{I}_n$ (as defined in \cite{HeningDang18c}). It is also interesting to note that by \cite[Proposition 3.1]{HeningDang18c} $\widetilde{\delta}(n)$ and $\mathcal{I}_n$ have the same sign and as only the positivity of these two parameters is important, working with one or the other doesn't matter.
\end{remark}

\subsection{Main results}


With the help of the previous parts, we can now state the stochastic persistence theorem.

\bthm \label{theorem_stoch_pers}
Suppose that $\sigma_1>0$ or $\sigma_n>0$ and that $\widetilde{\delta} (n) > 0$. Then
\bit
\ita $\PP_{inv}(\RR_{++}^n) = \{ \Pi\}$.

\itb $\Pi$ is absolutely continuous with respect to the Lebesgue measure of $\RR^n$ and its density is $\mathcal{C}^\infty$.

\itc For all $f \in L^1( \Pi)$ and $x \in \RR_{++}^n$, $\lim\limits_{t \to \infty} \Pi_t^x f = \Pi f$ a.s. and $(P_t)$ converges to $\Pi$ in total variation.

\eit
\ethm

The next result specifies the rate of convergence of $(P_t)$ toward $\Pi$.

\begin{theorem}\label{thm_rate_of_convergence}
Suppose that $\sigma_1>0$ or $\sigma_n>0$ and that $\widetilde{\delta} (n) > 0$.
\begin{itemize}
\item[(a)] If $a_{ii} = 0$ for some $i=2,\ldots,n$, then there exists $q_0>1$ such that for all $1<q<q_0$ there exists a proper $\CCC^2$ function $W_q : \RR_{++}^n \to \RR_+$ such that for all $ 1 \leq \beta \leq q$
$$\lim_{t\to\infty} t^{\beta-1} \left\|P_t(x,\cdot)-\Pi(\cdot)\right\|_{W_{\beta,q}}=0,\quad x\in\RR_{++}^n$$
where $W_{\beta,q} = W_q^{1-\beta/q}$.

\item[(b)] If $a_{ii} >0$ for all $i$, then there exist $\varsigma>0$ and a continuous proper function $\widehat{W} : \RR_{++}^n \to [1,\infty[$ such that
$$\lim_{t\to\infty} e^{\varsigma t} \left\|P_t(x,\cdot)-\Pi(\cdot)\right\|_{\widehat{W}}=0, \quad x\in\RR_{++}^n.$$
\end{itemize}
\end{theorem}

\begin{remark}
Theorem \ref{thm_rate_of_convergence} clearly implies convergence in total variation (take $\beta = q$ in $(a)$ and note that $\widehat{W} \geq 1$ in $(b)$).
\end{remark}


We have seen in the previous theorems that positivity of parameter $\widetilde{\delta}(n)$ ensures the persistence of all species. If we want to characterise extinction of at least one of them, we have to suppose that there exists $1\leq i\leq n$ such that $\widetilde{\delta}(i) \leq 0$. Intuitively if species $i$ goes extinct, species $i+1$ must become extinct too since its only source of food is species $i$. Arguing in the same manner, we intuitively conclude that species $i+2,\ldots,n$ also must become extinct. This will be done in Lemma \ref{Properties_delta}.

A. Hening and D. Nguyen already proved it in \cite{HeningDang18c,HeningDang18b} in the degenerate and non-degenerate setting :

\begin{theorem}[A. Hening, D. Nguyen 2018]\label{Thm_HN}
Suppose that there exists $1 \leq j^* < n$ such that $\widetilde{\delta}(j^*) >0$ and $\widetilde{\delta}(j^*+1) \leq 0$.
\begin{enumerate}
\item[(a)] If $a_{ii} = 0$ for some $i=2,\ldots,n$, then $\lim\limits_{t\to\infty} \frac1t \int_0^t \E_x[X_k(s)] ds = 0$, for $k=j^*+1,\ldots,n$ and all $x \in \RR_{++}^n$. Moreover, for every $\varepsilon>0$, there exists a compact set $K_\varepsilon \subset \RR_{++}^{(j^*)}$ such. that for all $x \in \RR_{++}^n$,
$$\liminf_{t \to \infty} \frac{1}{t}\int_0^t \P_x \left[ (X_1(s),\ldots,X_{j^*}(s))\in K_\varepsilon \right] ds\geq 1-\varepsilon.$$

\item[(b)] If $a_{ii} > 0$ for all $i=1,\ldots,n$, then $X_{j^*+1},\ldots,X_n$ go extinct almost surely exponentially fast and for every $\varepsilon>0$, there exists a compact set $K_\varepsilon \subset \RR_{++}^{(j^*)}$ such. that for all $x \in \RR_{++}^n$,
$$\liminf_{t \to \infty} \P_x \left[ (X_1(t),\ldots,X_{j^*}(t))\in K_\varepsilon \right] \geq 1-\varepsilon.$$
Moreover, if $\sigma_1,\ldots,\sigma_{j^*}>0$, there exists a unique invariant probability measure $\Pi$ supported by $\RR_{++}^{(j^*)}$ such that for all $x \in \RR_{++}^n$,
$P_t(x,\cdot)$ converges weakly to $\Pi(\cdot)$, i.e. for all continuous bounded functions $f : \RR_+^n \to \RR$,
$$\lim_{t\to \infty} P_tf(x) =\Pi f(x).$$

\end{enumerate}
\end{theorem}

\begin{remark}
Statement (a) of Theorem \ref{Thm_HN} is proved in \cite{HeningDang18b} in the case where $a_{ii} = 0$ for all $i \geq 2$. The proof easily extends verbatim to the present situation.
\end{remark}

The next theorem complements the second assertion of Theorem \ref{Thm_HN}.

\begin{theorem} \label{Thm_ext_1}
Suppose that there exists $1 \leq j^* <n$ such that $\widetilde{\delta}(j^*)>0$ and $\widetilde{\delta}(j^*+1) \leq 0$ and moreover that $\sigma_1>0$ or $\sigma_{j^*}>0$. Then
\bit
\ita $\PP_{inv}(\RR_+^n) = \PP_{inv}(\RR_+^{(j^*)})$ and $\PP_{inv}(\RR_{++}^{(j^*)}) = \{\Pi\}$.

\itb If $a_{ii}>0$ for all $i=1\ldots,n$, $P_t(x,\cdot)$ converges weakly to $\Pi(\cdot)$ for all $x \in \RR_{++}^n$.



\eit

\end{theorem}
\begin{remark}
In \cite{HeningDang18b}, when except for the lowest tropic level, species have no intraspecific competition (that is, $a_{ii}=0, i=2,\ldots,n$), the authors cannot show that $X_{k}(t), k\geq j^*$ converges to $0$ exponentially fast with probability 1 if $\widetilde{\delta}(j^*)>0$ and $\widetilde{\delta}(j^*+1)\leq 0$ because it is not easy to show the tightness of the family of random occupation measures when $a_{ii}=0, i=2,\ldots,n$.

Using newly-developed coupling techniques, one can obtain the exponential convergence of $X_{k}(t), k\geq j^*$ with probability 1. This result is beyond the scope of this paper and it will be reported in \cite{BNNpreprint}.
\end{remark}


\section{Stochastic persistence }

\subsection{Some mathematical tools} \label{Sect_Math_tolls_stoch_pers}

In this section, we introduce some tools that will be needed for the proof of Theorem \ref{theorem_stoch_pers}.  We consider the general stochastic differential equation given by \eqref{model_normal} under the assumption that $F_i$ is $\mathcal{C}^\infty$ but not necessarily of the form \eqref{fucntions_F_i}.

Let $L$ and $\Gamma_L$ respectively denote the {\it formal generator} and the {\it formal carré du champ}  defined by
$$Lg(x) := \sum_{i=1}^n x_i F_i(x) \frac{\partial g}{\partial x_i} (x) + \frac{1}{2} \sum_{i=1}^n \sigma_i^2 x_i^2 \frac{\partial^2 g}{\partial x_i^2} (x)$$
and
$$\Gamma_L (g)(x) := \sum_{i=1}^n \sigma_i^2 x_i^2 \left( \frac{\partial g}{\partial x_i} (x) \right)^2 $$
for $g : \RR_+^n \to \RR$ a $\mathcal{C}^2$ map (see \cite[Section 3.2]{B2018} for more details).

Let $C_b(\RR_+^n)$ be the set of real-valued bounded continuous functions on $\RR_+^n$. We let $\mathcal{L}$ be the {\it infinitesimal generator} of $(P_t)_{t \ge 0}$ on $C_b(\RR_+^n)$ and $\mathcal{D}(\mathcal{L}) \subset C_b( \RR_+^n)$ be its {\it domain}. Following \cite{B2018}, recall that $\mathcal{D}(\mathcal{L})$ is defined as the set of $f \in C_b(\RR_+^n)$ for which
\begin{enumerate}
\item $ \mathcal{L}f(x) :=\lim\limits_{t \to 0} \frac{P_tf(x)-f(x)}{t}$ exists for all $x \in \RR_+^n$;
\item $\mathcal{L}f \in C_b(\RR_+^n)$;
\item $\sup\limits_{0<t \le 1} \frac{1}{t} \|P_t f -f \|_\infty < \infty$.
\end{enumerate}

For a map $f \in \mathcal{D}(\mathcal{L})$ such that $f^2 \in \mathcal{D}(\mathcal{L})$, we let
$$\Gamma f := \mathcal{L}(f^2) - 2f \mathcal{L}f.$$
denote the {\it carrée du champ} operator.


The next result is Proposition $3.1$ in \cite{B2018} but for the second part of assertion $2$ for which we provide the proof in the Appendix \ref{U_complet_Stoch_pers}.

\begin{proposition} \label{Prop_fct_U}
Suppose that there exists a $\CCC^2$ proper map $U : \RR_+^n \to [1,\infty[$ and constants $\alpha,\gamma > 0$, $\beta \geq 0$ such that
\begin{align} \label{Inequality_function_U_in_Proposition}
LU \leq - \alpha U + \beta
\end{align}
and
\begin{align} \label{In_U2}
\Gamma_L(U) \leq \gamma  U^2.
\end{align}
Then,
\begin{itemize}
\item[1)] For each $x \in \RR_+^n$, there exists a unique strong solution $(X^x(t))_{t \geq 0} \subset \RR_+^n$ to \eqref{model_normal} with initial condition $X^x(0) = x$ and $X^x(t)$ is continuous in $(t,x)$. In particular, the process is $C_b(\RR_+^n)$-Feller, that is for each continuous bounded function $f : \RR_+^n \to \RR$, the mapping $(t,x) \mapsto P_t f(x)$ is continuous.

\item[2)] $\sup\limits_{t \geq 0} \E (U(X^x(t))) \leq \frac{\beta}{\alpha}$ for all $x \in \RR_+^n$ and  $\mu U = \int U d \mu \leq \frac{\beta}{\alpha}$ for all $\mu \in \PP_{inv}( \RR_+^n)$.

\item[3)] If $\CCC_c^2(\RR_+^n)$ is the set of $\CCC^2$ real valued maps with compact support, then $\CCC_c^2 (\RR_+^n) \subset \mathcal{D}^2(\mathcal{L})$, the set of functions $f$ such that $f,f^2 \in \mathcal{D}(\mathcal{L})$  and for all $f \in \CCC_c^2(\RR_+^n)$
$$\mathcal{L}f(x) = Lf(x) \quad \text{ and } \quad \Gamma(f)(x) = \Gamma_L(f)(x).$$

\item[4)] The set $\partial \RR_+^n$ is invariant under $(P_t)_{t \geq 0}$.

\item[5)] For $c>1$, the process
$$M_t(x) := U^{\frac{1}{2}}(X^x(t))-U^{\frac{1}{2}}(x) - c\int_0^t  U^{\frac{1}{2}}(X^x(s))ds, \quad t \ge 0$$
is a martingale that satisfies the strong law of large numbers, that is
$$\lim_{t \to \infty} \frac{M_t(x)}{t} = 0$$
$\P_x$ almost surely for all $x \in \RR_+^n$. There also exists $K \geq 0$ such that
$$LU^{\frac{1}{2}} \leq - c\, U^{\frac{1}{2}} + K.$$

\end{itemize}
\end{proposition}


In the proof of persistence, accessibility also plays a key role.

\begin{definition}
A point $y \in \RR_+^n$ is accessible from $x \in \RR_+^n$ if for every neighbourhood $U$ of $y$, there exists $t \geq 0$ such that $P_t(x,U)>0$.

We denote by $\textgoth{A}\hspace{1pt}_x$ the set of points $y$ that are accessible from $x$ and for $D \subset \RR_+^n$, we let $\textgoth{A}\hspace{1pt}_D = \underset{x \in D}{\bigcap} \textgoth{A}\hspace{1pt}_x$ be the set of accessible points from $D$.
\end{definition}

Let
\begin{align}\label{def_A_i}
A^0 (x) = (x_1\widetilde{F}_1(x),\ldots,x_n \widetilde{F}_n(x)) \quad \text{ and } \quad A^j(x) = \sigma_j x_j e_j, \; \; j=1,\ldots, n
\end{align}
where $e_j$ stands for the $j$-th vector of the canonical basis of $\RR^n$. Now, we consider the deterministic control system associated to \eqref{model_stratonovich}
\begin{align} \label{equa_control}
\dot{y} = A^0(y) + \sum_{j=1}^n u^j A^j(y)
\end{align}
where the control function $u = (u^1, \ldots, u^n) : \RR_+ \to \RR^n$ is at least piecewise continuous. We let $y(u,x,\cdot)$ denote the maximal solution to \eqref{equa_control} starting from $x$ and with control function $u$.

The next result is Proposition 5.3 in \cite{B2018}.

\begin{proposition} \label{prop_control}
Let $x,z \in \RR_+^n$, then $z \in \textgoth{A}\hspace{1pt}_x$ if and only if for every neighbourhood $O$ of $z$, there exists a control sequence $u$ such that $y(u,x,t) \in O$ for some $t\geq 0$.
\end{proposition}


We also need the so-called Hörmander condition. Let
\begin{align*}
E^1 & := \{A^1, \ldots, A^n\}, \\
E^k & := E^{k-1} \cup \left\{ [V,W] \mid V \in E^{k-1}, W \in E^{k-1} \cup\{A^0\} \right\} \qquad k \ge 2,
\end{align*}
where $[\cdot,\cdot]$ denotes the Lie bracket operator and $A^1,\ldots,A^n$ are defined by \eqref{def_A_i}. Recall that given two smooth vector fields $V,W : \RR^n \to \RR^n$ and $x \in \RR^n$,
$$[V,W](x) := DW(x) V(x)-DV(x) W(x),$$
where $DV(x) = \left(\frac{\partial V_i}{\partial x_j}(x)\right)_{i,j}$ stands for the Jacobian matrix of $V$ at $x$.

For $x \in \RR^n$ and $k \geq 1$, we let $E^k(x) := \left\{ V(x) \mid V \in E^k \right\}$.

\begin{definition}
We say that $\eqref{model_stratonovich}$ satisfies the strong Hörmander condition at $ x \in \RR_+^n$ if there exists $k \in \NN^*$ such that
$$span(E^k (x)) = \RR^n.$$

\end{definition}


For $\mu \in \mathcal{P}_{erg}(\partial \RR_+^n)$ provided $\widetilde{F}_i  \in L^1(\mu)$ we let
$$\mu \widetilde{F}_i := \int \widetilde{F}_i d\mu$$
denote the {\it invasion rate} of species $i$ with respect to $\mu$. See \cite{B2018} chapter $5$ or \cite{HeningDang18b} for more details.

Theorem \ref{theorem_stoch_pers} will be deduced from the next theorem which is a reformulation of Corollary $5.4$ in \cite{B2018}.

\begin{theorem} \label{Criter_thm_stoc_persi}
Suppose the hypotheses of Proposition \ref{Prop_fct_U} are satisfied and assume moreover that
\begin{align} \label{Criter_integrability_tilde_F}
\widetilde{F}_i \in L^1(\mu)\quad \forall \mu \in \mathcal{P}_{inv}(\partial\RR_+^n)\text{ and }i =1,\ldots ,n,
\end{align}
and that there exist positive numbers $p_1, \ldots, p_n$ such that for all $\mu \in \mathcal{P}_{erg}(\partial\RR_+^n)$
\begin{align} \label{Criter_p_i mu tilde F_i}
\sum_{i=1}^n p_i \mu \widetilde{F}_i > 0.
\end{align}
If there exists $x^* \in \textgoth{A}\hspace{1pt}_{\RR_{++}^n} \cap \RR_{++}^n$ satisfying the strong Hörmander condition, then
\begin{itemize}
\ita $\PP_{inv}(\RR_{++}^n) = \{ \Pi\}$ and $\Pi$ is absolutely continuous with respect to the Lebesgue measure of $\RR^n$.

\itb For all $f \in L^1( \Pi)$ and $x \in \RR_{++}^n$, $\lim\limits_{t \to \infty} \Pi_t^x f = \Pi f$ a.s.  and $(P_t)$ converge to $\Pi$ in total variation.

\itc The density of $\Pi$ is $\CCC^\infty$.
\end{itemize}
%
%

\end{theorem}

\brem \label{U_control_tilde_Fi}
By part $(2)$ of Proposition \ref{Prop_fct_U}, a sufficient condition for \eqref{Criter_integrability_tilde_F} is
$$\sum_{i=1}^n |\widetilde{F}_i | \leq a U + b$$
for some constants $a,b >0$.
\erem


\subsection{Proof of Theorem \ref{theorem_stoch_pers}} \label{Section_proof_thm_stoch_pers}


To prove the theorem, we assume that $F_i$ is given by \eqref{fucntions_F_i}.
%
%
%
Let
\beq \label{function_U}
U(x) = 1 + \sum_{i = 1}^n  c_i x_i
\eeq
where $c_1 >0$, $c_i = c_1 \prod\limits_{k=2}^i \frac{a_{k-1,k}}{a_{k,k-1}}$ for $i = 2,\ldots,n$.

\begin{lemma} \label{U_proof}
The function $U$ and the parameters
\begin{align*}
\alpha = \min_{i=2,\ldots,n} \{a_{i0}\}, && \beta = \alpha + \sup\limits_{x_1 >0} \{c_1 x_1(a_{10} - a_{11} x_1) + c_1 \alpha x_1 \}, &&\gamma = \max_{i=1,\ldots,n}\{\sigma_i^2\}
\end{align*}
satisfy hypotheses \eqref{Inequality_function_U_in_Proposition} and \eqref{In_U2} of Proposition \ref{Prop_fct_U}. Moreover
$$\widetilde{F}_i \in L^1(\mu)$$
$\forall \mu \in \mathcal{P}_{inv}(\partial\RR_+^n)\text{ and }i =1,\ldots ,n.$
\end{lemma}

\prf The first part of the proof is easy and left to the reader (see \cite[Lemma 3.2]{HeningDang18b} for a detailed one). The second one follows from Remark \ref{U_control_tilde_Fi}.

{\hfill$\square$\bigbreak}

We now relate the sign of $\widetilde{\delta}(n)$ to the persistence condition \eqref{Criter_p_i mu tilde F_i}.


\begin{proposition} \label{delta_pos_iif_existe_p_i_st_sum_pi_mu_lmada_i_pos}
$\widetilde{\delta}(n) >0$ is equivalent to condition \eqref{Criter_p_i mu tilde F_i} of Theorem \ref{Criter_thm_stoc_persi}.
\end{proposition}

\prf First suppose that $\widetilde{\delta}(n) >0$.  By Lemma 4 in \cite{SBA11}, condition \eqref{Criter_p_i mu tilde F_i} is equivalent to
\begin{align*}
\max_{i=1, \ldots, n} \mu \widetilde{F}_i >0, \tag{$\star$}
\end{align*}
for every $\mu \in \PP_{inv}(\partial \RR_+^n)$. By Lemma 3.5 of \cite{HeningDang18b}, the first part of the proof of Theorem $1.1$ $(i)$ in \cite{HeningDang18c} and Remark \ref{delta_cst_kappa}, $(\star)$ is verified.

Remark that in \cite{HeningDang18c,HeningDang18b} the aforesaid proofs do not require the non-degeneracy of the process, so they are still true in this setting.

Conversely, suppose \eqref{Criter_p_i mu tilde F_i} is true. By Theorem $5.1$ of \cite{B2018} and Lemma \ref{U_proof}, there exists $\mu \in \PP_{inv}(\RR_{++}^n)$ such that $\mu \widetilde{F}_i = 0$ for all $i=1,\ldots,n$. Now, write explicitly the system $\mu \widetilde{F}_i  = 0$ :
\begin{align*}
\left\{
     \begin{array}{ll}
     0 = \widetilde{a}_{10} -\widetilde{a}_{11} \int x_1\, \mu (dx)- \widetilde{a}_{12} \int x_2 \,\mu (dx) 								& i=1,\\
     0 = -\widetilde{a}_{i0} + \widetilde{a}_{i,i-1}\int x_{i-1} \,\mu (dx) - \widetilde{a}_{ii}\, \int x_i \,\mu (dx) - \widetilde{a}_{i,i+1} \int x_{i+1} \,\mu (dx) 	& i\neq 1,n , \\
     0 = -\widetilde{a}_{n0} + \widetilde{a}_{n,n-1} \int x_{n-1} \,\mu (dx) -\widetilde{a}_{nn} \int x_n \,\mu (dx) 					& i=n.
     \end{array}
     \right.
\end{align*}

So $\left(\int x_1 \mu(dx), \ldots, \int x_n \mu(dx)\right)$ is the unique solution of the previous system and by definition of $\mu$, this solution is strictly positive. By Proposition \ref{Prop_delt_pos_iff_F_positive_Zero}, $\widetilde{\delta}(n) >0 $.

{\hfill$\square$\bigbreak}


The next proposition is Theorem $5.3.1$ of \cite{HS98} for the intra-specific case and Exercise $5.3.2$ of \cite{HS98} for the non-intra-specific one.

\begin{proposition} \label{Prop_p_star_accessible}
Assume the equivalent conditions of Proposition \ref{Prop_delt_pos_iff_F_positive_Zero} are satisfied, then for the deterministic system $\dot{X}_i = X_i \widetilde{F}_i(X)$, $i = 1 \ldots, n$ (i.e. system \eqref{model_stratonovich} without noise), $p^*$ (as defined in Proposition \ref{Prop_delt_pos_iff_F_positive_Zero}) is globally asymptotically stable i.e. for every initial condition $x_0 \in \RR_{++}^n$, the solution $X^{x_0}(t)$ has the property that $X^{x_0}(t) \xrightarrow[t \to \infty]{} p^*$.
\end{proposition}

\begin{corollary} \label{p_star_acc}
Under the equivalent conditions of Proposition \ref{Prop_delt_pos_iff_F_positive_Zero}, $p^* \in \textgoth{A}\hspace{1pt}_{\RR_{++}^n}$.
\end{corollary}

\prf We start by rewriting the control equation for \eqref{model_stratonovich}
$$\dot{y}_i = y_i \, \widetilde{F}_i (y) + u^i\, y_i  \, \sigma_i \quad i=1,\ldots,n$$
where $u=(u^1,\ldots,u^n) : \RR_+ \to \RR^n$ is at least piecewise continuous.

If we take $u \equiv 0$, by Proposition \ref{Prop_p_star_accessible} the solution $y(0,x,t) \to p^*$ for every $x \in \RR_{++}^n$. We conclude by using Proposition \ref{prop_control}.

{\hfill$\square$\bigbreak}


Now, we discuss Hörmander's condition.

\begin{proposition} \label{Hormander_every_point_Rn}
Assume that either $\sigma_1 >0$ or $\sigma_n>0$. Then the strong Hörmander condition for \eqref{model_stratonovich} is satisfied for every $x \in \RR_{++}^n$.
\end{proposition}

\prf Recall that for $x \in \RR_+^n$, $A^1(x) = \sigma_1 x_1 e_1$ and $A^0 (x) = (x_1\widetilde{F}_1(x),\ldots,x_n \widetilde{F}_n(x))$. We immediately get for $x \in \RR_+^n$ that $DA^1 (x)$ and $DA^0(x)$ are the $n \times n$ matrix defined as
$$DA^1(x)_{ij} = \left\{ \begin{array}{ll}
\sigma_1 & \text{if } (i,j) = (1,1), \\
0 & \text{if } (i,j) \neq (1,1),
\end{array} \right.$$
and as
$$DA^0(x) = \begin{pmatrix}
\widetilde{F}_1 (x) + x_1\frac{\partial \widetilde{F}_1}{\partial x_1} (x) & x_1 \frac{\partial \widetilde{F}_1}{\partial x_2} (x) & 0 & \cdots & 0 \\
x_2 \frac{\partial \widetilde{F}_2}{\partial x_1} (x) & \widetilde{F}_2 (x) + x_2 \frac{\partial \widetilde{F}_2}{\partial x_2}(x) & x_2 \frac{\partial \widetilde{F}_2}{\partial x_3}(x) & \ddots & \vdots \\
0 & x_3 \frac{\partial \widetilde{F}_3}{\partial x_2} (x) & \ddots & \ddots & 0 \\
\vdots & \ddots & \ddots & \ddots & x_{n-1} \frac{\partial \widetilde{F}_{n-1}}{\partial x_n} (x)\\
0 & \cdots & 0 & x_n \frac{\partial \widetilde{F}_n}{\partial x_{n-1}}(x) & \widetilde{F}_n (x) + x_n \frac{\partial \widetilde{F}_n}{\partial x_n} (x)
\end{pmatrix}.$$
We now define recursively the family of vector fields $b^k \in E^k$, $k \geq 1$ by
$$\left\{\begin{array}{ll}
b^1 = A^1, \\
b^{k+1} = [b^k,A^0].
\end{array} \right.$$
We have for $x \in \RR_{++}^n$ that $b^k$ has the form
$$
b^k_k(x) = \sigma_1 x_1 \prod_{i=2}^k x_i \, \widetilde{a}_{i,i-1} \quad \text{ and } \quad b^k_j (x) = 0 \text{ for all } j >k.
$$
The assumption $\sigma_1 >0$ makes $\{b^1(x), \ldots, b^n(x)\}$ a basis of $\RR^n$ for every $ x \in \RR_{++}^n$. So Hörmander's condition holds true for every $x \in \RR_{++}^n$.

The case $\sigma_n >0$ is very similar and left to the reader.

{\hfill$\square$\bigbreak}

By Lemma \ref{U_proof}, the function $U$ defined by \eqref{function_U} satisfies the hypotheses of Proposition \ref{Prop_fct_U} and furthermore \eqref{Criter_integrability_tilde_F} is satisfied. Since $\widetilde{\delta}(n) >0$ and by Proposition \ref{delta_pos_iif_existe_p_i_st_sum_pi_mu_lmada_i_pos}, hypothesis \eqref{Criter_p_i mu tilde F_i} holds true. Under conditions $\sigma_1 >0$ or $\sigma_n >0$ and $\widetilde{\delta}(n) >0$, Corollary \ref{p_star_acc} and Proposition \ref{Hormander_every_point_Rn} make the equilibrium $p^*$ an accessible point satisfying the strong Hörmander condition. Proposition \ref{Hormander_every_point_Rn} proves moreover that the strong Hörmander condition holds at every $x \in \RR_{++}^n$.

Then hypotheses of Theorem \ref{Criter_thm_stoc_persi} hold true and Theorem \ref{theorem_stoch_pers} follows.




\section{Rate of convergence}\label{sec_rate_of_cv}

\subsection{General main results on the rate of convergence}

As in Section \ref{Sect_Math_tolls_stoch_pers}, we consider the general stochastic differential equation given by \eqref{model_normal} under the assumption that $F_i$ is $\mathcal{C}^\infty$ but not necessarily of the form \eqref{fucntions_F_i}.

We suppose here that there exists a map $U : \RR_+^n \to \RR $ satisfying the hypotheses of Proposition \ref{Prop_fct_U} and the additional conditions :
\beq\label{exc0}
\liminf_{\|x\|\to\infty} \frac{U(x)}{\ln\|x\|}>0,
\eeq
\beq\label{exc1}\limsup_{\|x\|\to\infty}\left(LU(x)+ p_0\sum_{i=1}^n |F_i(x)|\right)<0 \text{ for some } p_0>0
\eeq
and
\beq\label{exc2} \sum_{i=1}^n |F_i(x)|\leq C U^{d_0}(x) \text{ for some } d_0\geq 1 \text{ and } C>0.
\eeq
We also assume that the persistence condition of Theorem \ref{theorem_stoch_pers} is satisfied, that is there exist $p_1,\dots, p_n>0$ such that
\beq\label{exc3}\sum_{i=1}^n p_i \mu\widetilde{F}_i>0 \text{ for every $\mu \in \PP_{erg}(\partial \RR_+^n)$.}
\eeq
Remark that we can assume that $\sum_{i=1}^n p_i$ is sufficiently small without loss of generality.

Consider the function
$$
V(x)= U(x)-\sum_{i=1}^n p_i \ln x_i.
$$
When  $\sum_{i=1}^n p_i$ is sufficiently small, we deduce from \eqref{exc0} and \eqref{exc1} that the function $V(x)$ satisfies
$$V(x)>0 \text{ for any } x\in \RR^n_{++} \,\text{ and }\lim_{\|x\|\to\infty} V(x)=\infty$$
%
and
\beq\label{e3.1}
LV(x)\leq h_1\Ind_{\{\|x\|\leq M\}} -h_2 \Ind_{\{\|x\|\geq M\}}, \quad x\in \RR^n_{++}
\eeq
%
for some positive constants $h_1, h_2, M$.

Let us define the constant $q_0 >1$ such that
\begin{align}\label{def_q0}
-\alpha +\frac{q_0-1}2\gamma=0
\end{align}
where $\alpha$ and $\gamma$ are as in \eqref{Inequality_function_U_in_Proposition} and \eqref{In_U2}.

Finally, for $q \in \left]1, \min\{q_0, \frac{q_0+2}2\}\right[$, let
\begin{align}\label{fct_W}
W_q = \left\{ \begin{array}{ll}
V^q + CU^q & \text{if } 1<q \leq 2, \\
V^q+C U^{2q-2} & \text{if } q>2, \\
\end{array}\right.
\end{align}
and
\begin{align}\label{fct_hat_W}
\widehat{W}(x) = \left(\frac{U(x)}{\prod_{i=1}^n x_i^{p_i}}\right)^{\eps^*}
\end{align}
for some $C,\varepsilon^* >0$.
Remark that $W_q >0$ and that $\widehat{W} \geq 1$ when $\sum_{i=1}^n p_i$ is sufficiently small.


We can now state two new general theorems for the rates of convergence.


\begin{theorem}\label{Thm_polynomial_rate}
Assume that there exists a map $U : \RR_+^n \to \RR $ satisfying the hypotheses of Proposition \ref{Prop_fct_U}, \eqref{exc0}, \eqref{exc1}, \eqref{exc2} and \eqref{exc3}.
Suppose moreover that there exists $x^* \in \textgoth{A}\hspace{1pt}_{\RR_{++}^n} \cap \RR_{++}^n$ which satisfies the strong Hörmander condition. Then for all $q \in \left]1, \min\{q_0, \frac{q_0+2}2\}\right[$ and for all $ 1 \leq \beta \leq q$,
$$\lim_{t\to\infty} t^{\beta-1} \left\|P_t(x,\cdot)-\Pi(\cdot)\right\|_{W_{\beta,q}}=0,\quad x\in\RR_{++}^n$$
where $W_{\beta,q} = W_q^{1-\beta/q}$ and $W_q$ is defined by \eqref{fct_W}.
\end{theorem}


\begin{theorem}\label{Thm_exponential_rate}
Assume that there exists a map $U : \RR_+^n \to \RR $ satisfying hypotheses of Proposition \ref{Prop_fct_U}, \eqref{exc2} and \eqref{exc3}. Suppose that \eqref{exc0} and \eqref{exc1} are strengthened to
\beq\label{exc0_strong}
\liminf_{\|x\|\to\infty} \frac{\ln U(x)}{\ln\|x\|}>0
\eeq
and
\begin{align}\label{cv_exp_hyp_strong}
\limsup_{\|x\|\to\infty}\left(L \ln U(x)+ p_0\sum_{i = 1}^n |F_i(x)|\right)<0 \text{ for some } p_0>0
\end{align}
respectively.
Assume moreover that there exists $x^* \in \textgoth{A}\hspace{1pt}_{\RR_{++}^n} \cap \RR_{++}^n$ which satisfies the strong Hörmander condition. Then there exists $\varsigma>0$ such that
$$\lim_{t\to\infty} e^{\varsigma t} \left\|P_t(x,\cdot)-\Pi(\cdot)\right\|_{2\widehat{W}}=0, \quad x\in\RR_{++}^n$$
where $\widehat{W}$ is defined by \eqref{fct_hat_W}.
\end{theorem}

Before proving these theorems, we will apply them in the next part to prove the polynomial or exponential rate of convergence in the Lotka-Volterra food chain model.

\subsection{Proof of Theorem \ref{thm_rate_of_convergence}}

To prove Theorem \ref{thm_rate_of_convergence}, we will use Theorem \ref{Thm_polynomial_rate} and \ref{Thm_exponential_rate}.

\begin{itemize}
\item[(a)] Suppose that $a_{ii} =0$ for some $i=2,\ldots,n$. By Lemma \ref{U_proof}, $U$ satisfies the hypotheses of Proposition \ref{Prop_fct_U}. We have from \eqref{Inequality_function_U_in_Proposition} that
$$LU(x) + p_0 \sum_{i=1}^n |F_i(x)| \leq -\alpha U + \beta + p_0 \sum_{i=1}^n |F_i(x)|.$$
Since $F_i$ and $U$ are linear functions, then for $p_0$ sufficiently small \eqref{exc1} holds true and \eqref{exc0} is trivially verified. We also deduce that for each $d_0 \geq 1$, there exists a constant $C>0$ such that \eqref{exc2} holds. Since $\widetilde{\delta}(n) >0$, by Proposition \ref{delta_pos_iif_existe_p_i_st_sum_pi_mu_lmada_i_pos}, \eqref{exc3} is immediately verified.
Moreover by Corollary \ref{p_star_acc} and Proposition \ref{Hormander_every_point_Rn}, $p^*\in \textgoth{A}\hspace{1pt}_{\RR_{++}^n} \cap \RR_{++}^n$ satisfies the strong Hörmander condition. 

Thus, by Theorem \ref{Thm_polynomial_rate}, for all $q \in \left]1, \min\{q_0, \frac{q_0+2}2\}\right[$ and for all $ 1 \leq \beta \leq q$
$$\lim_{t\to\infty} t^{\beta-1} \left\|P_t(x,\cdot)-\Pi(\cdot)\right\|_{W_{\beta,q}}=0,\quad x\in\RR_{++}^n$$
where $q_0$ is defined by \eqref{def_q0}, $W_q$ by \eqref{fct_W} and $W_{\beta,q} = W_q^{1-\beta/q}$.

\item[(b)] Suppose that $a_{ii} >0$ for all $i$. The verification of the hypotheses of Proposition \ref{Prop_fct_U} and of inequalities \eqref{exc2} and \eqref{exc3} is basically the same as in the case $(a)$. \eqref{exc0_strong} follows from the form of $U$. To prove \eqref{cv_exp_hyp_strong}, first remark that
\begin{align*}
L \ln U(x) & = \frac{LU(x)}{U(x)} - \frac{1}{2} \frac{\Gamma_L(U)(x)}{U^2(x)} \\
& = \frac{c_1 x_1 (a_{10} - a_{11}x_1) - \sum_{i=2}^n c_i a_{i0} x_i - \sum_{i=2}^n c_i a_{ii} x_i^2}{1+ \sum_{i=1}^n c_i x_i} - \frac{1}{2} \frac{\sum_{i=1}^n c_i^2\sigma_i^2 x_i^2}{\left(1+\sum_{i=1}^n c_i x_i\right)^2}.
\end{align*}
It thus becomes clear that $L\ln U(x)$ is bounded above by a linear function of $\|x\|$ with a negative slope. Then \eqref{cv_exp_hyp_strong} holds true for $p_0>0$ sufficiently small.
Since $\widetilde{\delta}(n)>0$, by Corollary \ref{p_star_acc} and Proposition \ref{Hormander_every_point_Rn}, $p^*\in \textgoth{A}\hspace{1pt}_{\RR_{++}^n} \cap \RR_{++}^n$ satisfies the strong Hörmander condition.

Thus, by Theorem \ref{Thm_exponential_rate}, there exists $\varsigma>0$ such that
$$\lim_{t\to\infty} e^{\varsigma t} \left\|P_t(x,\cdot)-\Pi(\cdot)\right\|_{2\widehat{W}}=0, \quad x\in\RR_{++}^n$$
where $\widehat{W}$ is defined by \eqref{fct_hat_W}.

\end{itemize}

This concludes the proof.

\begin{remark}
In the case $(a)$, i.e. $a_{ii}=0$ for some $i=2,\ldots,n$, if $\sigma_1>0$ and $\sigma_i=0$, for $i=2$ $\ldots,n$, then we can easily verify that for any $q>0$,
$$L U^q(x)\leq \beta_q-\alpha_q U^q(x)$$ for some $\beta_q,\alpha_q>0$.
Thus, in view of Proposition \ref{lm3.3}, the degree of $t^{\beta-1}$ in Theorem \ref{thm_rate_of_convergence} part $(a)$ can be any positive number.

\end{remark}

\brem
The positivity of $a_{ii}$ allows us to choose a stronger Lyapunov function. It implies that $\left(X_t^x\right)$ enter in a compact exponentially fast and this is why the semi-group $\left(P_t\right)$ converges to the stationary distribution exponentially fast.
\erem


\subsection{Proofs of Theorem \ref{Thm_polynomial_rate} and Theorem \ref{Thm_exponential_rate}}\label{prf_thm_cv_general}

We start by Theorem \ref{Thm_polynomial_rate}.
In order to obtain the polynomial convergence rate, we will show that for each $q \in \left]1, \min\{q_0, \frac{q_0+2}2\}\right[$, there is a $T>0$ such that the proper function $W_q:\RR_{++}^n\to\RR_+$ defined by \eqref{fct_W} satisfies
\begin{align}\label{inequality_W_q}
P_T(W_q(x)) -W_q(x) \leq- \kappa W_q^{(q-1)/q}(x)+K \Ind_{\{W_q \leq R\}}(x)
\end{align}
for some $\kappa,K,R>0$. We will also show that
the Markov chain $\left(X^x(kT)\right)_{k\in\NN}$ is aperiodic, irreducible and all compact set in $\RR^n_{++}$ is petite. Then by \cite[Theorem 3.6]{jarner2002polynomial}, for all $ 1 \leq \beta \leq q$
$$\lim_{t\to\infty} t^{\beta-1} \left\|P_t(x,\cdot)-\Pi(\cdot)\right\|_{W_{\beta,q}}=0,\quad x\in\RR_{++}^n$$
where $W_{\beta,q} = W_q^{1-\beta/q}$.

To proceed, let us recall some
concepts and results needed to prove the main theorem.
Let ${\bf\Phi}=(\Phi_0,\Phi_1,\dots)$ be a discrete-time Markov chain on a general state space $(E,\mathcal{E})$, where $\mathcal{E}$ is a countably generated $\sigma$-algebra.
Denote by $\mathcal{P}$ the Markov transition kernel for ${\bf\Phi}$.
If there is a non-trivial $\sigma$-finite positive measure $\varphi$ on $(E,\mathcal{E})$ such that for
any $A\in\mathcal{E}$ satisfying $\varphi(A)>0$ we have
$$\sum_{k=1}^\infty \mathcal{P}^k(x, A)>0,\quad x\in E$$
where  $ \mathcal{P}^k$ is the $k$-step transition kernel of ${\bf\Phi}$, then the Markov chain ${\bf\Phi}$ is called $\varphi$-\textit{irreducible}.
It is known that if ${\bf\Phi}$ is $\varphi$-irreducible, then there exist a positive integer $d$ and disjoint subsets $E_0,\dots,
E_{d-1} \subset E$ such that for all $i=0,\dots, d-1$ and all $x\in E_i$, we have
$$\mathcal{P}(x,E_j)=1 \text{ where } j=i+1 \text{ (mod $d$)}$$
and
$$ \varphi \left(E\setminus \bigcup_{i=0}^{d-1}E_i\right)=0.$$
The smallest positive integer $d$ satisfying the above condition is called the {\it period} of ${\bf\Phi}.$
An \textit{aperiodic} Markov chain is a chain with period $d=1$.

A set $C\in\mathcal{E}$ is called \textit{petite}, if there exist a non-negative sequence $(a_k)_{k\in\NN}$ with $\sum_{k=1}^\infty a_k=1$
and a nontrivial positive measure $\nu$ on $(E,\mathcal{E})$
such that
$$\sum_{k=1}^\infty a_k \mathcal{P}^k(x, A)\geq\nu(A),\quad x\in C, A\in\mathcal{E}.$$


\begin{lemma}\label{acc_Horm_discret}
	If there exists $p^* \in \textgoth{A}\hspace{1pt}_{\RR_{++}^n} \cap \RR_{++}^n$ satisfying the strong Hörmander condition, then for any $T>0$, the Markov chain $\left(X^x(kT)\right)_{k\in\NN}$ is aperiodic, irreducible and all compact sets in $\RR^n_{++}$ are petite.
\end{lemma}

\prf
	Due to the strong Hörmander condition, by \cite{ichihara1974classification} the transition probability $P_t( x,\cdot)$ of $\left(X^x(t)\right)$ has a smooth density, denoted by $p(t,x,y)$. Let $q^* \in \RR_{++}^n$ and $T>0$ such that $p(T, p^*, q^*)>0$.
	Since $p(t,x,y)$ is smooth, we can find $r>0$ such that
	\begin{equation}\label{edd0}
	\inf\left\{p(T,x,y)\mid x\in B(p^*, r), y\in B(q^*, r)\right\}=:\delta_r>0,
	\end{equation}
	where $B(x,r)$ is the open ball with radius $r$ centered at $x$.
	
	
	By accessibility of $p^*$, for each compact set $K\subset\RR_{++}^n$ there exists $k_K>0$ such that $$P_{kT}(x, B(p^*, r/2))>0$$ for any $x\in K$, $k\geq k_K$.
	Due to the compactness of $K$ and the Feller property of $\left(X^x(t)\right)$, we can easily obtain that
	\begin{equation}\label{edd2}
	\inf_{x\in K} \{P_{kT}(x, B(p^*, r))\}:=\varsigma_{k, K}>0, \text{ for any }k\geq k_K.
	\end{equation}
	Applying the Chapman–Kolmogorov equation and using \eqref{edd0} and \eqref{edd2}, we have
	\begin{equation}\label{edd3}
	\inf_{x\in K} \{P_{(k+1)T}(x, C)\}\geq \varsigma_{k, K}\delta_r m(C\cap B(q^*,r)), \quad k\geq k_K
	\end{equation}
 for any measurable subset $C\subset\RR_{++}^n$, where $m(\cdot)$ is the Lebesgue measure on $\RR_{++}^n$.
 Obviously, \eqref{edd3} implies that the Markov chain $\left((X^x(kT)\right)_{k\in\NN}$ is irreducible and every compact set $K$ is petite.
 Moreover, $\left((X^x(kT)\right)_{k\in\NN}$ is aperiodic because if it is not, then there is $d>1$ and disjoint subsets $E_0,\dots,
 E_{d-1}$ such that for all $i=0,\dots, d-1$ and $ x\in E_{i}$, $$P_T(x, E_j)=1,\quad j=i+1 \text{ (mod $d$)}.$$
 This implies that for any $k \in  \NN$ and $x\in E_{i}$
 	\begin{equation}\label{edd4}
 	P_{(kd+1)T}(x, E_j)=1,\quad  j=i+1 \text{ (mod $d$)}.
 	\end{equation}
 Now, observe that at least one of the two sets $E_0\cap B(q^*,r)$ and $B(q^*,r)\setminus E_0$ has nonzero Lebesgue measure.
 If $m(E_0\cap B(q^*,r))>0$, for any $x\in E_0$, letting $K$ be a compact set containing $x$, then we have from \eqref{edd3} that
 $$P_{(kd+1)T}(x, E_0)\geq  \varsigma_{k, K}\delta_rm(E_0\cap B(q^*,r))>0, \quad k\geq k_K$$
 which contradict the fact that $E_0, E_1$ are disjoint and $P_{(kd+1)T}( x, E_1)=1$.
 Likewise, if $m(B(q^*,r)\setminus E_0)>0$, for any $x\in E_{d-1}$, letting $K$ be a compact set containing $x$, we have from \eqref{edd3} that
 $$ P_{(kd+1)T}( x, B(q^*,r)\setminus E_0)\geq \varsigma_{k, K}\delta_rm(B(q^*,r)\setminus E_0)>0, \quad k\geq k_K$$
 which contradict $P_{(kd+1)T}(x, E_0)=1$ and $x\in E_{d-1}$.
 The contradiction argument above shows that the Markov chain $\left((X^x(kT)\right)_{k\in\NN}$ is aperiodic.

{\hfill$\square$\bigbreak}

Now, we want to prove that for all $q \in \left]1, \min\{q_0, \frac{q_0+2}2\}\right[$ $W_q$ satisfies \eqref{inequality_W_q}. The proof will follow the approach in \cite{BNNpreprint}. To proceed, we first need the next technical results. We start by next lemma that provided very useful inequalities.

\begin{lemma}\label{lm3.0}
	Let $p>0$. There exists $c_p>0$ such that for any $a>0$ and $x\in\RR$ we have
	\beq\label{lm3.0-e1}
	|a+x|^{1+p}\leq \begin{cases}
		a^{1+p}+(1+p)a^{p}x+c_p|x|^{1+p} &\text{ if } p\leq 1,\\
		a^{1+p}+(1+p)a^{p}x+a^p+c_p\left(|x|^{2p}+1 \right) & \text{ if } p>1.
	\end{cases}
	\eeq
	As an application,
	for a random variable $Y$ and a constant $k>0$, one has that there exists $\tilde K_1 >0$ such that
	\beq\label{lm3.0-e2}
	\E |k+Y|^{1+p}\leq
	\begin{cases}
		k^{1+p}+(1+p)k^p{\E Y}+\tilde K_1\E |Y|^{1+p} &\text{ if } p\leq 1,\\
		k^{1+p}+(1+p)k^p{\E Y}+k^p+\tilde K_1\left(\E |Y|^{2p}+1 \right) &\text{ if } p>1.
	\end{cases}
	\eeq
	Moreover,  there exists $c_{p, b}>0$ depending only on $p,b>0$ such that if $x + a \geq 0$ then
	\beq\label{lm3.0-e3}
	(a+x)^{1+p}-b(a+x)^p \leq  \begin{cases}
		a^{1+p}+ (1+p)a^p x-\frac{b}2 a^p+c_{b,p}(|x|^{1+p}+1) &\text{ if } p\leq 1,\\
		a^{1+p}+ (1+p)a^p x-\frac{b}2 a^p+a^p+c_{b,p}(|x|^{2p}+1)&\text{ if } p>1.
	\end{cases}
	\eeq
	
\end{lemma}

\prf
	If $x>a/2$ then 	
	$$(a+x)^{1+p}-a^{1+p}-(1+p)a^px< c_p x^{1+p}.$$
	If $x<-a/2$ then
	$$|a+x|^{1+p}-a^{1+p}-(1+p)a^px\leq c_p|x|^{1+p}.$$
	If $p\leq 1$ we are done for this two cases. For $p>1$, since $a>0$, it is enough to remark that $c_p |x|^{1+p} \leq a^p + c_p( |x|^{2p} +1)$.
	\\ If $|x| \leq a/2$, consider the function
	$$f(x)=(a+x)^{1+p}-a^{1+p}-(1+p)a^px.$$
	We have $f(0)=0, f'(0)=0, f''(x)=(1+p)p(a+x)^{p-1}$.
	By Taylor's expansion, there exists $0\leq \xi_x\leq 1$ such that
	$$
	f(x)= (1+p)p\frac{x^2(a+\xi_x x)^{p-1}}{2}
	$$
	If $0<p\leq 1$ then
	$x^2(a+\xi_x x)^{p-1}\leq x^2 (2|x|)^{p-1}\leq c_p |x|^{p+1}$.
	If $p>1$ then
	we have from Young's inequality that
	$$x^2(a+\xi_x x)^{p-1}\leq x^2 (3/2a)^{p-1}\leq a^p + c_p |x|^{2p}.$$
	Thus, \eqref{lm3.0-e1} has been proved. \eqref{lm3.0-e2} straightforwardly follows.
	To prove \eqref{lm3.0-e3}, we consider two cases.
	If $0\leq a+x\leq 2^{-1/p} a$ then
	we have $(a+x)^{1+p}\leq 2^{-(1+p)/p} a^{1+p}$ and the inequality is true by choosing a constant $c_{b,p}$ sufficiently large that
	$$2^{-(1+p)/p} a^{1+p}+ \frac{b}2 a^p\leq a^{1+p}+c_{b,p}.
	$$
	If $a+x\geq 2^{-1/p} a$ then $-b(a+x)^p\leq -\frac{b}2 a^p$
	and \eqref{lm3.0-e3} follows \eqref{lm3.0-e1}.
	
	{\hfill$\square$\bigbreak}

Let
$2\rho:=\inf\{\sum_{i=1}^n p_i \mu\widetilde{F}_i: \mu\in\PP_{inv}(\partial \RR_+^n)\} >0$
and let $n^*$ be a positive  integer such that
$ (n^*-1)h_2-2h_1\geq \frac\rho2$.

\begin{lemma}\label{lm3.1}
Suppose that $U$ satisfies hypotheses of Proposition \ref{Prop_fct_U} and inequalities \eqref{exc0}, \eqref{exc1}, \eqref{exc2} and \eqref{exc3}.
 	Then there exist $T^*>0, \delta>0$ such that
	$$
	\E_x\int_0^T LV(X(s))ds \leq -\rho T $$
	for any  $T\in [T^*, n^*T^*]$, $x\in\RR_{++}^n$, $\|x\|\leq M$, and $\dist(x,\partial\RR_+^n)\leq\delta$.
	
\end{lemma}
\prf
	The lemma is well known; see  \cite[Proposition 8.2]{B2018} or \cite[Proposition 4.1]{HN18}.
	
	{\hfill$\square$\bigbreak}
	
\begin{lemma}\label{lm3.2}
	Suppose $U$ satisfies the hypotheses of Proposition \ref{Prop_fct_U}. Let $q_0$ like in \eqref{def_q0}, i.e. $q_0>1$ satisfies $-\alpha +\frac{q_0-1}2\gamma=0$ where $\alpha$ and $\gamma$ are as in \eqref{Inequality_function_U_in_Proposition} and \eqref{In_U2}.
	For any $1<q<q_0$,
		\beq\label{e2-lm3.2}L U^q(x)\leq k_{1q} - k_{2q} U^q(x), \quad x\in\RR_+^n
	\eeq
	for some positive constants $k_{1q}, k_{2q}$.
	As a result,
	\beq\label{e1-lm3.2}
	\E_x U^q(X(t))\leq \frac{k_{1q}}{k_{2q}} + U^q(x) e^{-k_{2q}t}, \quad t\geq 0.
	\eeq
	On the other hand,
	\beq\label{e3-lm3.2}
	\E_x U^{p}(X(t))\leq e^{K_pt} U^{p}(x)
	\eeq
	for any $p\geq 1$ and some constant $K_p >0$.
\end{lemma}

\prf
	The estimate \eqref{e2-lm3.2} is straight forward calculations but we give here the proof for completeness. For $1<q<q_0$, by \eqref{Inequality_function_U_in_Proposition} and \eqref{In_U2} we have
	\begin{align*}
	LU^q(x) &= q U^{q-1} (x)\, LU(x) + \frac12	q(q-1) \, U^{q-2} (x) \, \Gamma_L (U)(x) \\
	&\leq  q U^q (x)\, \left(-\alpha +\frac{q-1}{2} \gamma\right) +q \beta  U^{q-1}(x).
	\end{align*}
	Since $q<q_0$, then there exist $k_{1q},k_{2q} >0$ such that
	\begin{align*}
	L U^q(x)\leq k_{1q} - k_{2q} U^q(x), \quad x\in\RR_+^n.
	\end{align*}

	Then \eqref{e1-lm3.2} follows by a standard application of It\^o's formula to the function $e^{k_{2q}t} U^q(X(t))$.
	Finally, it is easy to check that
	$$
	LU^{p}(x)\leq K_p U^{p}(x)
	$$
	for some constant $K_p>0$. An application of It\^o's formula and Gronwall's inequality derives \eqref{e3-lm3.2}. A more detailed proof can be found in \cite[Section $9.3$]{B2018} for \eqref{e1-lm3.2} and \eqref{e3-lm3.2}.
	
{\hfill$\square$\bigbreak}
\begin{proposition}\label{lm3.3}
	Assume $U$ satisfies the hypotheses of Proposition \ref{Prop_fct_U} and inequalities \eqref{exc0}, \eqref{exc1}, \eqref{exc2} and \eqref{exc3}.
	Let $q$ be any number in the interval $\left]1, \min\{q_0, \frac{q_0+2}2\}\right[$.
	There is $\kappa^*>0$ and $C, C^{*}>0$ such that
	$$
	\E_x [CU^q(X(n^*T^*))+V^q(X(n^*T^*))]\leq CU^q(x)+V^q(x)-\kappa^* [CU^q(x)+V^q(x)]^{\frac{q-1}q} + C^{*}
	$$
	if $1<q\leq 2$
	and
		$$
	\E_x [CU^{2q-2}(X(n^*T^*))+V^q(X(n^*T^*))]\leq CU^{2q-2}(x)+V^q(x)-\kappa^* [CU^{2q-2}(x)+V^q(x)]^{\frac{q-1}q} + C^{*}
	$$
	if $q>2$.

\end{proposition}

\prf
First we assume that $1<q\leq 2$.
	In the sequel, $C^*$ is a generic constant depending on $T^*, M, n^*$ but independent of $x\in\RR_{++}^n$.
	$C^*$ can differ from line to line.
	Suppose $X(0)=x$. We have from It\^o's formula that
	$$
	V(X(t))=V(x)+\int_0^t LV(X(s))ds+M(t).
	$$
	It can be seen that $M(t)$ is a martingale with quadratic variation:
	\beq\label{e3-lm3.3}\langle M(t)\rangle=\int_0^t \Gamma_L(V)(X(s))ds\leq K \int_0^t U^{2}(X(s))ds
	\eeq
	for some constant $K$ independent of $t$.
%
	
	Note that by \eqref{exc2}, $|LV(x)|\leq CU^{d_0}(x)$ which together with \eqref{e3-lm3.3} and \eqref{e3-lm3.2} and applications of It\^o's isometry and Hölder's inequality implies that
	\beq\label{e4-lm3.3}
	\E_x\left|\int_0^t LV(X(s))ds\right|^q\leq C^* U(x)^{qd_0} \,\,\text{ and }\,\E_x\left|M(t)\right|^q\leq C^* U^{q}(x), \quad t\leq n^*T^*.
	\eeq
	We have from \eqref{e4-lm3.3} and \eqref{lm3.0-e2} that
	\beq\label{e0-lm3.3}
	\begin{aligned}
		\E_x [V(X(t))]^q\leq&  V^q(x) + q\left[\E_x \int_0^t LV(X(s))ds\right] V^{q-1}(x)+ C^*\E_x\left|\int_0^t LV(X(s))ds+M(t)\right|^q\\
		\leq &  V^q(x) + q\left[\E_x \int_0^t LV(X(s))ds \right] V^{q-1}(x)+ C^*(1+U^{qd_0}(x)), \quad t\leq n^*T^*.
	\end{aligned}
	\eeq
	Thus, if $\|x\|\leq M $ and $ \dist(x,\partial \RR_+^n)\leq \delta$, we have  $\E_x \int_0^t LV(X(s))ds\leq -\rho t$, $t\in[T^*, n^*T^*]$. As a result
	\beq\label{e5-lm3.3}
	\begin{aligned}
		\E_x [V(X(T))]^q
		\leq &  V^q(x) - q\rho T V^{q-1}(x)+ C^*(1+U^{qd_0}(x)), \quad T\in[T^*, n^*T^*].
	\end{aligned}
	\eeq
	On the other hand, we have from $$V(X(T))=V(x)+\int_0^T LV(X(s))ds+M(T)\leq V(x)+h_1 T+M(T)$$  and \eqref{lm3.0-e2} that
	\beq\label{e6-lm3.3}
	\begin{aligned}
		\E_x [V(X(T))]^q
		\leq &  V^q(x) + qh_1 T V^{q-1}(x)+ C_2(1+U^{qd_0}(x)), \quad T\leq n^*T^*.
	\end{aligned}
	\eeq
	for some $C_2$ independent of $x$.
	Since $V(x)$ is bounded on the set $\{x\in\RR_{++}^n: \|x\|\leq M, \dist(x,\partial \RR_{++}^n)\geq\delta\}$, by
	combining \eqref{e5-lm3.3} and \eqref{e6-lm3.3} for $\|x\| \leq M$, we obtain that
	\beq\label{e7-lm3.3}
	\begin{aligned}
		\E_x [V(X(T))]^q
		\leq &  V^q(x) - q\rho T V^{q-1}(x)+ C^*,\quad T\in[T^*, n^*T^*].
	\end{aligned}
	\eeq

	
	Define 	 $$\xi=\inf\{t\geq 0: \|X(t)\|\leq M\}\wedge (n^*T^*).$$
	From now on, we suppose that $\|x\| \geq M$. For $t\leq \xi$, we have
	\beq\label{e9-lm3.3}
	V(X(t))= V(x)+\int_0^tLV(X(s))ds +M(t)\leq V(x)- h_2t+M(t).
	\eeq
	We have from \eqref{e7-lm3.3} and the strong Markov property of $X(t)$ that
	\beq\label{e1-lm3.3}
	\begin{aligned}
		\E_x& \left[\Ind_{\{\xi \leq T^*(n^*-1)\}} V^q(X(n^*T^*))\right]\\[0.5ex]
		\leq &\E_x \left[\Ind_{\{\xi \leq T^*(n^*-1)\}} \left[V^q(X(\xi))+C^*\right] \right]\\[0.5ex]
		&-\E_x\left[\Ind_{\{\xi \leq T^*(n^*-1)\}}q\rho(n^*T^*-\xi)V^{q-1}(X(\xi))\right]\\[0.5ex]
		\leq& \E_x \left[\Ind_{\{\xi \leq T^*(n^*-1)\}} (V(x)+M(\xi))^q+C^* \right] \\[0.5ex]
		&-q\rho T^*\E_x\left[\Ind_{\{\xi \leq T^*(n^*-1)\}}(V(x)+M(\xi))^{q-1}\right]\\[0.5ex]
		&\text{ (due to \eqref{e9-lm3.3})}\\
		\leq& \E_x\left[\Ind_{\{\xi \leq T^*(n^*-1)\}}\left( V^q(x)-\frac{q\rho T^*}2 V^{q-1}(x)+ qM(\xi)V^{q-1}(x)+C^*(|M(\xi)|^q+1)\right)\right]\\[0.5ex]
		&\text{ (\eqref{lm3.0-e3} is applied here).}
	\end{aligned}
	\eeq
	If $T^*(n^*-1)\leq \xi \leq T^*n^*$ we have
	\begingroup
	\allowdisplaybreaks
	\begin{align*}
		\E_x &\left[\Ind_{\{\xi \geq T^*(n^*-1)\}} V^q(X(n^*T^*))\right]\\[0.5ex]
		\leq &\E_x \left[\Ind_{\{\xi \geq T^*(n^*-1)\}} V^q(X(\xi))+C^*\right] \\[0.5ex]
		&+qh_1\E_x\left[\Ind_{\{\xi \geq T^*(n^*-1)\}}(n^*T^*-\xi)V^{q-1}(X(\xi))\right] \\[0.5ex]
		&\text{ (due to \eqref{e6-lm3.3} and the strong Markov property)}\\[0.5ex]
		\leq& \E_x \left[\Ind_{\{\xi \geq T^*(n^*-1)\}} [(V(x)+M(\xi)- h_2\xi)^q+C^*]\right] \\[0.5ex]
		&+qh_1 T^*\E_x\left[\Ind_{\{\xi \geq T^*(n^*-1)\}}(V(x)+M(\xi)- h_2\xi)^{q-1}\right] \\[0.5ex]
		&\text{ (because of \eqref{e9-lm3.3})}\\[0.5ex]
		\leq& \E_x\left[\Ind_{\{\xi \geq T^*(n^*-1)\}}\left( V^q(x)-q h_2\xi V^{q-1}(x)+qM(\xi) V^{q-1}(x) + C^*\left(|M(\xi)|+1\right)^q\right)\right]\\[0.5ex]
		&+ 2^{q-1}qh_1 T^*\E_x\left[\Ind_{\{\xi \geq T^*(n^*-1)\}} \left( V^{q-1}(x)+|M(\xi)|^{q-1}\right)\right]\\[0.5ex]
		& \text{ (\eqref{lm3.0-e1} and the inequality $|x+y|^{q-1}\leq 2^{q-1}(|x|^{q-1}+|y|^{q-1})$ are applied here)} \\[0.5ex]
		\leq & \E_x \left[ \Ind_{\{\xi \geq T^*(n^*-1)\}} \left( V^q(x)-\frac{q\rho T^*}2 V^{q-1}(x)+q M(\xi)V^{q-1}(x) +C^*\left(|M(\xi)|+1\right)^q \right)\right] \\
		& \text{ (since } (n^*-1)h_2-2^{q-1}h_1\geq \frac\rho2 \text{ and } \xi\geq T^*(n^*-1)).\\[-2ex]
		& \stepcounter{equation}\tag{\theequation}\label{e2-lm3.3}
	\end{align*}
	\endgroup
	As a result, by adding \eqref{e1-lm3.3} and \eqref{e2-lm3.3} and noting that $\E_x M(\xi)=0$, we have
	\beq\label{e10-lm3.3}
	\begin{aligned}
		\E_x V^q(X(n^*T^*))\leq&  V^q(x)-q\frac\rho2 T^* V^{q-1}(x)+ C^* \E_x (|M(\xi)|+1)^q\\
		\leq& V^q(x)-q\frac\rho2 T^* V^{q-1}(x) + C^* U^{q}(x)
	\end{aligned}
	\eeq
	where the inequality $\E_x (|M(\xi)|+1)^q\leq C^* U^{q}(x)$ comes from an application of
	the Burkholder-Davis-Gundy Inequality, Hölder's inequality and \eqref{e3-lm3.3} and \eqref{e3-lm3.2}.
	
	From \eqref{e1-lm3.2}, we have
	\beq\label{e11-lm3.3}
	\begin{aligned}
		\E_x U^q(X(n^*T^*))\leq&  U^q(x)-\left(1-e^{-k_{2q} n^*T^*}\right) U^q(x) + \frac{k_{1q}}{k_{2q}}.
	\end{aligned}
	\eeq
	Combining \eqref{e10-lm3.3} and \eqref{e11-lm3.3}, we can easily get that
	\beq\label{e12-lm3.3}
	\E_x \left[V^q(X(n^*T^*))+CU^q(X(n^*T^*))\right]
	\leq V^q(x)+C U^q(x)-\kappa^* [V^q(x)+CU^q(x)]^{(q-1)/q} + C^{*},
	\eeq
	for some $ \kappa^*>0, C^{*}>0$ and sufficiently large $C$.

If $q>2$ then $q_0>2$ and $2q-2<q_0$ since $q<\min\{q_0,\frac{q_0+2}2\}$.
We carry the proof in the same manner with use of the inequalities \eqref{lm3.0-e1}, \eqref{lm3.0-e2} and \eqref{lm3.0-e3} for the case $p=q-1>1$ and we obtain that
	\beq\label{e13-lm3.3}
\begin{aligned}
	\E_x V^q(X(n^*T^*))\leq&  V^q(x)-q\frac\rho2 T^* V^{q-1}(x)+ V^{q-1}(x)+ C^* \E_x (|M(\xi)|+1)^{2q-2}\\
	\leq& V^q(x)-(q\frac\rho2 T^*-1) V^{q-1}(x) + C^* U^{2q-2}(x)\\
	\leq&V^q(x)-(q\frac\rho4 T^*) V^{q-1}(x) + C^* U^{2q-2}(x)
\end{aligned}
\eeq
when we choose $T^*$ sufficiently large. (Note that $\rho$ and $q$ do not depend on the choice of $T^*$ so we can choose $T^*>\frac{2}{q\rho}$.)
Since $2q-2<q_0$, we have from \eqref{e1-lm3.2} that
\beq\label{e11-lm3.3_v2}
\begin{aligned}
	\E_x U^{2q-2}(X(n^*T^*))\leq&  U^{2q-2}(x)-\left(1-e^{-k_{2,(2q-2)} n^*T^*}\right) U^{2q-2}(x) + \frac{k_{1,(2q-2)}}{k_{2,(2q-2)}}.
\end{aligned}
\eeq
Combining two displays above, we obtain
\beq\label{e12-lm3.3_v2}
\begin{aligned}
\E_x \left[V^q(X(n^*T^*))+CU^{2q-2}(X(n^*T^*))\right]&\\
 &\hspace{-4cm}\leq V^q(x)+C U^{2q-2}(x)
 -\kappa^* [V^q(x)+CU^{2q-2}(x)]^{(q-1)/q} + C^{*},
\end{aligned}
\eeq
for some $ \kappa^*>0, C^{*}>0$ and sufficiently large $C$.

{\hfill$\square$\bigbreak}

Fix $ q \in \left]1, \min\{q_0, \frac{q_0+2}2\}\right[$, then $W_q$ defined by \eqref{fct_W} is proper and for $R>0$, the set $\{W_q\leq R\}$ is compact and by Lemma \ref{acc_Horm_discret} also petite. Then for $0< \hat\kappa < \kappa^*$, by Proposition \ref{lm3.3} there exists $R>0$ such that
$$P_{n^*T}W_q - W_q \leq -\hat\kappa W_q^{(q-1)/q} + C^* \Ind_{\{ W_q\leq R\}}. $$
Using Lemma \ref{acc_Horm_discret} and Theorem $3.6$ of \cite{jarner2002polynomial}, we get for all $ 1 \leq \beta \leq q$
$$\lim_{t\to\infty} t^{\beta-1} \left\|P_t(x,\cdot)-\Pi(\cdot)\right\|_{W_{\beta,q}}=0,\quad x\in\RR_{++}^n$$
where $W_{\beta,q} = W_q^{1-\beta/q}$. This concludes the proof of Theorem \ref{Thm_polynomial_rate}.

\bigskip
\bigskip
\bigskip
We pass now to the proof of Theorem \ref{Thm_exponential_rate} and we first need the next proposition.

\begin{proposition}
	Suppose that hypotheses of Proposition \ref{Prop_fct_U}, \eqref{exc2}, \eqref{exc3}, \eqref{exc0_strong}, \eqref{cv_exp_hyp_strong} hold.
	Then
	there is $\kappa^*\in(0,1)$ and $K^*>0$, $\eps^*>0$ such that
	$$
	\E_x \left(\frac{U(X(t))}{\prod_{i=1}^n X_i^{p_i}(t)}\right)^{\eps^*}\leq C^*+\kappa^*\left(\frac{U(x)}{\prod_{i=1}^n x_i^{p_i}}\right)^{\eps^*}, \quad x\in \RR_{++}^n.
	$$
\end{proposition}


 \prf
	The result is basically proved in \cite{B2018, HN18}.
	
{\hfill$\square$\bigbreak}

In other words, this Proposition shows that
$$P_t \widehat{W}(x) \leq K^* + \kappa^* \widehat{W}(x), \quad x\in \RR_{++}^n$$
for some constants $K^*>0$, $\kappa \in (0,1)$ independent of $x$ and $t$ and for $\widehat{W}$ defined by \eqref{fct_hat_W}. For $R$ large enough and since $\widehat{W}$ is proper, we have
$$P_t \widehat{W}(x) \leq \hat\kappa^* \widehat{W}(x) +  \hat K^* \Ind_{\{ \widetilde{W}\leq R\}}(x), \quad x\in \RR_{++}^n$$
for some $\hat\kappa^* \in (\kappa,1)$ and a $ \hat K^*>0$. Then by Lemma \ref{acc_Horm_discret}, the set $\{ \hat W\leq R\}$ is petite. By the second corollary of Theorem $6.2$ of \cite{meyn1992stability}, there exists $\varsigma>0$ such that
$$\lim_{t\to\infty} e^{\varsigma t} \left\|P_t(x,\cdot)-\Pi(\cdot)\right\|_{2\widehat{W}}=0, \quad x\in\RR_{++}^n.$$
This concludes the proof of Theorem \ref{Thm_exponential_rate}.


\section{Extinction}


In Section \ref{main_result}, we have seen that intuitively, if species $j^*$ goes extinct then species $j^*+1,\ldots,n$ do, too.
Next lemma confirms that this intuition is correct. Remark that A. Hening and D.Nguyen already proved it in \cite{HeningDang18c,HeningDang18b} with other tools since they didn't have the constant $\widetilde{\delta}(n)$.


\begin{lemma} \label{Properties_delta}
\begin{itemize}
\item[(a)] If $\widetilde{\delta}(n) > 0$, then $\widetilde{\delta}(n-1) >0$.
\item[(b)] If $\widetilde{\delta}(n-1) \leq 0$, then $\widetilde{\delta}(n) < 0$.
\item[(c)] If $\widetilde{\delta}(n) \leq 0$, then there is no invariant probability measure supported by $\RR_{++}^n$.
\end{itemize}
\end{lemma}

\prf
First, notice that $\widetilde{\delta}(n)$ can be rewritten as
\begin{align*}
\widetilde{\delta}(n) &= \widetilde{a}_{10} \prod_{j=2}^n \widetilde{a}_{j,j-1} - \sum_{k=2}^n \widetilde{a}_{k0} \prod_{l=k+1}^n \widetilde{a}_{l,l-1} \sum_{\alpha \in A_1^{k-1}} \prod_{j=1}^{k-1} \widetilde{a}_{j,\alpha(j)}\\
&= \widetilde{a}_{n,n-1}\left[ \widetilde{a}_{10} \prod_{j=2}^{n-1} \widetilde{a}_{j,j-1} - \sum_{k=2}^{n-1} \widetilde{a}_{k0} \prod_{l=k+1}^{n-1} \widetilde{a}_{l,l-1} \sum_{\alpha \in A_1^{k-1}} \prod_{j=1}^{k-1} \widetilde{a}_{j,\alpha(j)}\right]
- \widetilde{a}_{n0} \sum_{\alpha \in A_1^{n-1}} \prod_{j=1}^{n-1} \widetilde{a}_{j,\alpha(j)} \\
&= \widetilde{a}_{n,n-1}\widetilde{\delta}(n-1)- \widetilde{a}_{n0} \sum_{\alpha \in A_1^{n-1}} \prod_{j=1}^{n-1} \widetilde{a}_{j,\alpha(j)}.
\end{align*}
Then parts $(a)$ and $(b)$ become clear. Suppose there exists $\mu \in \PP_{inv}(\RR_{++}^n)$, then by same arguments as in the proof of Proposition \ref{delta_pos_iif_existe_p_i_st_sum_pi_mu_lmada_i_pos}, $\widetilde{\delta}(n) >0$ and $(c)$ follows.

{\hfill$\square$\bigbreak}

\begin{remark} \label{Remark_properties_delta}
The same conclusions holds true with $n$ replaced by $i$ in last lemma and in particular if $\widetilde{\delta}(i) \leq 0$ there is no invariant probability measure supported by $\RR_{++}^{(i)}$.
In particular, this implies that if there exists $1 \leq i \leq n$ such that $\widetilde{\delta}(i) \leq 0$, then $\widetilde{\delta}(k)<0$ for all $k=i+1, \ldots,n$ and furthermore, there is no invariant probability measures supported by $\RR_{++}^{(k)}$ for all $k=i,\ldots,n$.
\end{remark}

\subsection{Proof of Theorem \ref{Thm_ext_1}}

Since $\widetilde{\delta}(j^*)>0$ and $\sigma_1>0$ or $\sigma_{j^*}>0$, by Theorem \ref{theorem_stoch_pers} we obtain that
\begin{enumerate}
\item[$(I)$] $\PP_{inv}(\RR_{++}^{j^*}) = \{ \Pi\}$.

\item[$(II)$] 
$(P^{(j^*)}_t)$ converge to $\Pi$ in total variation where  $(P^{(j^*)}_t)$ stands for the transition kernel restricted to the first $j^*$ species.
\end{enumerate}
Point $(a)$ then follows from Remark \ref{Remark_properties_delta} and $(I)$.

By the form of the Lyapunov function $\widehat{W}$ constructed in Theorem \ref{Thm_exponential_rate} and by the second Corollary of Theorem $6.2$ of \cite{meyn1992stability}, $(P_t)$ converge weakly uniformly on each compact set of $\RR_{++}^{(j^*)}$, i.e. for each continuous bounded function $f : \RR_{++}^{(j^*)} \to \RR_+$ and each compact set $K\subset\RR_{++}^{(j^*)}$,
$$\lim_{t\to \infty}\left(\sup_{x\in K} \left|P_tf(x)-\Pi f(x)\right|\right) = 0 .$$
By \cite{HeningDang18c}, we know that species $x_1,\ldots,x_{j^*}$ are persistent in probability, that is for any $\varepsilon>0$, there exists a compact set $K_\varepsilon\subset\RR_{++}^{(j^*)}$ such that for all $x \in \RR_{++}^n$,
$$\liminf_{t \to \infty} \P_x \left[ (X_1(t),\ldots,X_{j^*}(t))\in K_\varepsilon \right] \geq 1-\varepsilon.$$
Moreover, always by \cite{HeningDang18c}, species $x_{j^*+1},\ldots,x_n$ goes extinct almost surely. Then by the latter, the Feller property and by mimicking the proof of Theorem 1.1 $(iii)$ of \cite{HeningDang18c}, we get point $(b)$, that is $P_t(x,\cdot)$ converges weakly to $\Pi(\cdot)$ for all $x \in \RR_{++}^n$.

{\hfill$\square$\bigbreak}


\section{Appendix}


\subsection{Proof of Proposition \ref{Prop_fct_U} point 2}\label{U_complet_Stoch_pers}

The first part is Proposition 3.1 in \cite{B2018} and the second one is Lemma $3.2$ of \cite{HeningDang18b}, but we give the proof for completeness.

Let $\tau_k = \inf \left\{ t \geq 0 \mid U(X^x(t)) \geq k\right\} $ for $x$ fixed. We write $a\wedge b$ for the infimum between $a$ and $b$. By Dynkin's formula and hypothesis \eqref{Inequality_function_U_in_Proposition},
\begin{align*}
\E_x \left[ U( X(\tau_k \wedge t) )\right]  \; & = \; U(x) + \E_x \left[ \int_0^{\tau_k \wedge t} LU(X(s)) ds \right] \\
& \leq \; U(x) + \beta \, \E_x [\tau_k \wedge t] - \alpha \, \E_x \left[ \int_0^{\tau_k \wedge t} U(X(s)) ds \right] .
\end{align*}
By letting $k \to \infty$, we obtain
$$\E_x \left[ U( X(t) )\right] \leq U(x) + \beta t - \alpha \, \E_x \left[ \int_0^{t} U(X(s)) ds \right].$$
By definition of a Markov kernel and the Fubini-Tonelli theorem, we finally have
$$P_t U(x) \leq U(x) +\beta t - \alpha \int_0^t P_s U(x) ds.$$
By the same arguments as in the proof of Theorem 2.2 in \cite{B2018}, we get
$$P_t U(x) \leq  e^{-\alpha t} \left( U(x) - \frac{\beta}{\alpha} \right) + \frac{\beta}{\alpha}.$$
This implies that for all $ x \in \RR_+^n$
$$\limsup\limits_{t\to \infty} P_t U(x) \leq \frac{\beta}{\alpha}.$$

We pass now to the proof of the second part. Let $\mu \in \PP_{inv}(\RR_+^n)$ and $K >0$. By Fatou's lemma
\begin{align*}
\mu \left( K \wedge U \right) \; & = \;  \mu \left[ P_t ( K \wedge U) \right] \\
& = \;\lim\limits_{t \to \infty} \mu \left[ P_t ( K \wedge U) \right] \\
& \leq \; \mu \left[ \limsup\limits_{t \to \infty} P_t ( K \wedge U) \right] \\
& \leq \; \frac{\beta}{\alpha}.
\end{align*}
By letting $K \to \infty$, we obtain that for every $\mu \in \PP_{inv}(\RR_+^n)$, $\mu U \leq \frac{\beta}{\alpha}$.

{\hfill$\square$\bigbreak}

\subsection{Proof of Proposition \ref{Prop_delt_pos_iff_F_positive_Zero}} \label{sec_proof_delta_pos}

For $i \in \{1,\ldots, n-2\}$ and $\alpha \in A_{i+1}^n$, $\beta \in A_{i+2}^n$, we define
$$\widehat{\alpha}(j) = \left\{ \begin{array}{ll}
\alpha(j) & \text{ if } j \neq i, \\
i  & \text{else,}
\end{array} \right. \qquad \text{ and } \qquad \widehat{\beta} = \beta \, (i \; i+1). $$
Then, the map $\widehat{\color{white}\alpha} : A_{i+1}^n \cup A_{i+2}^n \to A_i^n$ is a bijection.

We now write explicitly the system $\widetilde{F}(x)=0$ :
$$ \left\{
     \begin{array}{rcll}
      \widetilde{a}_{10} & = & \widetilde{a}_{11} \, x_1 + \widetilde{a}_{12} \, x_2, & \\
      \widetilde{a}_{i,i-1} \, x_{i-1}  & = & \widetilde{a}_{i0} + \widetilde{a}_{ii} \, x_i +  \widetilde{a}_{i,i+1} \,x_{i+1} &  i= 2, \ldots,n-1, \\
      \widetilde{a}_{n,n-1} \, x_{n-1} & = & \widetilde{a}_{n0} + \widetilde{a}_{nn} \,  x_n. &
     \end{array}
     \right.
$$

One can check by induction that
\begin{align*}
x_i  = \prod_{j=i+1}^n \frac{1}{\widetilde{a}_{j,j-1}} &\left[  x_n \sum_{\alpha \in A_{i+1}^n} \prod_{j=i+1}^n \widetilde{a}_{j,\alpha(j)}  + \widetilde{a}_{i+1,0} \prod_{j=i+2}^n \widetilde{a}_{j,j-1} \right. \\
&   +  \left.\sum_{k = i+2}^n \widetilde{a}_{k0} \prod_{l= k+1}^n \widetilde{a}_{l,l-1} \sum_{\alpha \in A_{i+1}^{k-1}} \prod_{j=i+1}^{k-1} \widetilde{a}_{j,\alpha(j)} \right]
\end{align*}
for $i = 1, \ldots, n-1$.

Now, we determine the value of $x_n$. By the first equation of the system and the previous relation, we have
\begingroup
\allowdisplaybreaks
\begin{align*}
\widetilde{a}_{10} \prod_{j=2}^n \widetilde{a}_{j,j-1} & = x_1 \, \widetilde{a}_{11} \prod_{j=2}^n \widetilde{a}_{j,j-1} + x_2 \, \widetilde{a}_{12} \, \widetilde{a}_{21} \prod_{j=3}^n \widetilde{a}_{j,j-1} \\
%
%
& = x_n \left[ \sum_{\alpha \in A_{2}^n} \widetilde{a}_{11} \prod_{j=2}^n \widetilde{a}_{j,\alpha(j)} + \sum_{\beta \in A_{3}^n} \widetilde{a}_{12} \, \widetilde{a}_{21} \prod_{j=3}^n \widetilde{a}_{j,\beta(j)} \right] \\
%
%
& +  \widetilde{a}_{20}\, \widetilde{a}_{11} \prod_{j=3}^n \widetilde{a}_{j,j-1}   +  \sum_{k = 3}^n \widetilde{a}_{k0} \prod_{l= k+1}^n \widetilde{a}_{l,l-1} \sum_{\alpha \in A_2^{k-1}} \widetilde{a}_{11} \prod_{j=2}^{k-1} \widetilde{a}_{j,\alpha(j)} \\
%
%
& +  \widetilde{a}_{30}\, \widetilde{a}_{12}\, \widetilde{a}_{21} \prod_{j=4}^n \widetilde{a}_{j,j-1}   +  \sum_{k = 4}^n \widetilde{a}_{k0} \prod_{l= k+1}^n \widetilde{a}_{l,l-1} \sum_{\beta \in A_3^{k-1}} \widetilde{a}_{12} \, \widetilde{a}_{21} \prod_{j=3}^{k-1} \widetilde{a}_{j,\beta(j)} \\
%
%
& = x_n \sum_{\alpha \in A_1^n} \prod_{j=1}^n \widetilde{a}_{j,\alpha(j)} + \sum_{k = 2}^n \widetilde{a}_{k0} \prod_{l= k+1}^n \widetilde{a}_{l,l-1} \sum_{\alpha \in A_1^{k-1}}  \prod_{j=1}^{k-1} \widetilde{a}_{j,\alpha(j)}.
\end{align*}
\endgroup
The last equality follows from the $\widehat{\color{white}\alpha}$ bijection. Then we have
\begin{align*}
 x_n \sum_{\alpha \in A_1^n} \prod_{j=1}^n \widetilde{a}_{j,\alpha(j)} = \widetilde{\delta}(n).
\end{align*}
As the sum to the left is strictly positive, we have that $x_n > 0$ iff $\widetilde{\delta} (n) >0$. We conclude the proof by noting that $x_n >0$ iff $x >0$.

\section*{Acknowledgments}
Michel Benaim and Antoine Bourquin are supported in part by the SNF grant 200020-196999.
Dang H. Nguyen is supported in part by NSF through the grant DMS-1853467.

We thank two anonymous referees for their useful comments and valuable suggestions.


\bibliographystyle{amsplain}
\bibliographystyle{imsart-nameyear}
\bibliographystyle{nonumber}
\bibliography{Degenerate_LV_Food_Chain_Biblio}
\end{document}